\documentclass{article}

\usepackage[english]{babel}

\usepackage[letterpaper,top=2cm,bottom=2cm,left=3cm,right=3cm,marginparwidth=1.75cm]{geometry}

\usepackage{amsfonts}
\usepackage{graphicx}
\usepackage[all]{xy}
\usepackage{amsmath,textcomp,amssymb,geometry,graphicx,enumerate}
\usepackage{mathrsfs}
\usepackage[colorlinks=true, allcolors=blue]{hyperref}
\title{On a specific problem of partition of graph}
\author{Peisheng Yu}

\begin{document}
\maketitle

\begin{abstract}
In this short article, we consider a problem about $2$-partition of the vertices of a graph. If a graph admits such a partition into some 'small' graphs, then the number of edges cross an arbitrary cut of the graph $e(S,S^{c})$ has a nice lower bound. 
\end{abstract}

\section{The problem}

To state the problem, we need first to introduce a special type of graph. For a positive real number $c$, define $c$-small graph as follow: if $M$ is the adjacent matrix of the graph, then such a graph is called $c$-small if $cJ-M$ is non-negatively definited, where $J$ is the all-one matrix. In other words, the inequality $x^{t}Mx \leq c \cdot x^{t} J x = c \cdot (\sum_{i} x_{i})^{2}$ always holds for any real vector $x$. 
\\ \hspace{\fill} \\
The problem we are going to consider in this article is as follow:
Given a graph $G=(V,E)$, suppose there is a $2$-partition $\mathcal{A}$ of the vertices $V$, such that:
\\(a) For every set $A \in \mathcal{A}$, we have $|A| \geq 2$.
\\(b) For every $ \{ u, v \} \subset V$, there is a unique $A \in \mathcal{A}$ such that $ \{ u, v \} \subset A$.
\\If for all $A \in \mathcal{A}$, the induced subgraph $G[A]$ is a $c$-small graph, then the following property holds: for any subset $S \subset V$, we have:
\begin{equation}
    e(S,S^{c}) \geq \lambda (c) \cdot min(e(S), e(S^{c}))
\end{equation}
In (1), $S^{c}=V \setminus S$ denotes the complement of $S$ in $V$, and $e(S,S^{c})$ denotes the number of edges between the two sets $S$ and $S^{c}$. $\lambda (c)$ is a positive real number depends only on $c$. In fact, as we will see, $\lambda (c)$ can take to be $\frac{2(1-c)}{1+c}$.

\section{The Proof}

The idea follows that of M. Fiedler \cite{fiedler1973algebraic}. It is well known that the second smallest eigenvalue of the Laplacian matrix (also called Fiedler value) of a graph detects its algebraic connectivity \cite{fiedler1989laplacian}. However, in most cases, the precise value of Fielder value is hard to obtain. To solve the problem, we estimate $x^{t}Lx$ and bound it from the below by some other quantities. 
\\ \hspace{\fill} \\
Let $A=(a_{ij})$ be the adjacent matrix of $G$, and $L=diag(d_{1},d_{2},\cdots, d_{n})-A=D-A$ be the Laplacian matrix of $G$. Here we assume $|V|=n$, and $d_{i}=deg(v_{i})$ for $v_{i} \in V$. For a real vector $x$, consider $x^{t}Lx$: we have $x^{t}Lx=\sum_{i,j}a_{ij}(x_{i}-x_{j})^{2}$. If $|S|=pn$, and $|S^{c}|=qn$, for some positive real number $p,q$, with $p+q=1$. Next, let:
\begin{equation}
    x_{v}=
    \begin{cases}
    q, & v \in S \\
    -p, & v \in S^{c}
\end{cases}
\end{equation}
Then, $a_{ij}(x_{i}-x_{j})^{2}=0$ if $v_{i}$ and $v_{j}$ belong to the same set $S$ or $S^{c}$; $a_{ij}(x_{i}-x_{j})^{2}=1$ if $v_{i}$ and $v_{j}$ belong to different sets. Thus, substitute above $x=(x_{v})_{v \in V}$ into $x^{t}Lx$, we obtain:
\begin{equation}
    e(S,S^{c})=x^{t}Lx
\end{equation}
On the other side, by direct computation, we have:
\begin{equation}
    x^{t}Lx=\sum_{k=1}^{n}d_{k}x_{k}^{2}-x^{t}Ax
\end{equation}
Since there is a $2$-partition for $V$, each edge $e\in E$ is belong to the edge set of the subgraph $G[A]$ for a unique $A \in \mathcal{A}$. If we write the partition $\mathcal{A}=\{ A_{1}, A_{2}, \cdots, A_{m} \}$, and denote the adjacent matrix of those $\{ G[A_{1}], G[A_{2}], \cdots, G[A_{m}] \}$ by $M_{1}, M_{2}, \cdots, M_{m}$, then we have $A= M_{1}+M_{2}+ \cdots + M_{m}$. By our assumption, each $G[A_{i}]$ is $c$-small graph, thus $x^{t}M_{i}x_{i} \leq c\cdot (\sum_{v\in A_{i}}x_{v})^{2}$. Summing up, we get:
\begin{gather}
    x^{t}Ax \leq c \cdot \sum_{i=1}^{m}(\sum_{v \in A_{i}}x_{v})^{2} = \nonumber \\
    c \cdot \sum_{i=1}^{m}\sum_{u,v \in A_{i}}x_{u}x_{v} \leq c \cdot \sum_{k=1}^{n}d_{k}x_{k}^{2}+2c \sum_{1 \leq i <j \leq n}x_{i}x_{j} 
\end{gather}
Combining this with (4), we get:
\begin{equation}
    x^{t}Lx \geq (1-c) \cdot \sum_{k=1}^{n}d_{k}x_{k}^{2}-2c \sum_{1 \leq i < j \leq n}x_{i}x_{j}
\end{equation}
Using (2), we compute that:
\begin{gather}
    \sum_{1 \leq i < j \leq n}x_{i}x_{j}=\binom{pn}{2}q^{2} + \binom{qn}{2}p^{2} - p^{2}q^{2}n^{2}= \nonumber \\
    \frac{pn(pn-1)}{2}q^{2} + \frac{qn}{2}p^{2} - p^{2}q^{2}n^{2}= -\frac{1}{2}pqn
\end{gather}
and also:
\begin{gather}
    \sum_{k=1}^{n}d_{k}x_{k}^{2}=\sum_{v\in S}d_{v}x_{v}^{2} + \sum_{v \in S^{c}}d_{v}x_{v}^{2} \nonumber \\
    = q^{2} \sum_{v \in S}d_{v} + p^{2} \sum_{u \in S^{c}}d_{u}
\end{gather}
Note that we also have:
\begin{gather}
    \sum_{v \in S}d_{v}=2e(S)+e(S,S^{c}) \\
    \sum_{u \in S^{c}} d_{u}=2e(S^{c})+e(S,S^{c})
\end{gather}
Combining (8), (9) and (10), we obtain:
\begin{equation}
    \sum_{k=1}^{n}d_{k}x_{k}^{2}=2q^{2}e(S)+2p^{2}e(S^{c})+(p^{2}+q^{2})e(S,S^{c})
\end{equation}
Putting (3), (6), (7) and (11) together:
\begin{equation}
    e(S,S^{c}) \geq (1-c)(p^{2}+q^{2})e(S,S^{c})+2(1-c)q^{2}e(S)+2(1-c)p^{2}e(S^{c})+cpqn
\end{equation}
(12) is equivalent to:
\begin{equation}
    (1-(1-c)(p^{2}+q^{2}))e(S,S^{c}) \geq 2(1-c)\cdot(q^{2}e(S)+p^{2}e(S^{c}))+cpqn
\end{equation}
The last term $cpqn$ is positive, omitting it:
\begin{equation}
    e(S,S^{c}) \geq (1-c)\cdot \frac{2p^{2}+2q^{2}}{1-(1-c)(p^{2}+q^{2})}\cdot min(e(S),e(S^{c}))
\end{equation}
Since $p+q=1$, $p^{2}+q^{2} \geq \frac{1}{2}$. We further get:
\begin{equation}
    e(S,S^{c}) \geq \frac{2(1-c)}{1+c} \cdot min(e(S),e(S^{c}))
\end{equation}

\section{Some Discussion}

\subsection{Examples of c-small graphs}
We give some examples $c$-small graphs in this section.
\\ \hspace{\fill} \\
Star graph, $c=\frac{1}{2}$. This is because $2x_{1}(x_{2}+ \cdots + x_{k}) \leq \frac{1}{2}(x_{1}+x_{2}+\cdots + x_{k})^{2}$. Thus, star graphs are $\frac{1}{2}$-small.
\\ \hspace{\fill} \\
Complete bipartite graph, $c=\frac{1}{2}$. The reason is the same as above: if $I \cup J =V$ is a partition, then $2(\sum_{v \in I}x_{v})(\sum_{u \in J}x_{u})\leq \frac{1}{2}(\sum_{v\in V}x_{v})^{2}$.
\\ \hspace{\fill} \\
Complete $s$-partite graph, $c=\frac{s-1}{s}$. This is because $2\sum x_{I}x_{J} \leq \frac{s-1}{s}(\sum x_{I})^{2}$ (On the left side there are $s$-many such $x_{I}$, with $x_{I} = \sum_{i\in I} x_{i}$).

\subsection{A refinement}
Since we omit the term $cpqn$ from (13) to (14), we can actually get a tighter bound if taking this term into consideration. Denote $e=min(e(S),e(S^{c}))$ to ease the notation. Then from (13):
\begin{equation}
    e(S,S^{c})\geq \frac{2(1-c)(1-2p+2p^{2})e+c(p-p^{2})n}{c+2(1-c)(p-p^{2})}
\end{equation}
If we denote the right side of (16) by $f(p)$, in order to get rid of $p$, we need to obtain the minimum value of $f(p)$:
\begin{equation}
    \frac{df(p)}{dp}=\frac{2p-1}{c+2(1-c)(p-p^{2})}\cdot (4(1-c)e-c^{2}n)
\end{equation}
We separate into two cases: 
\\ \hspace{\fill} \\
Case 1: $e > \frac{c^{2}}{4(1-c)}n$, then $f(p)$ reaches its minimum value at $p=\frac{1}{2}$. Therefore:
\begin{equation}
    e(S,S^{c}) \geq \frac{2(1-c)}{1+c}\cdot e + \frac{c}{4}n
\end{equation}
\\Case 2: $e \leq \frac{c^{2}}{4(1-c)}n$, then $f(p)$ reaches its minimum at $p=0$ or $p=1$. Therefore:
\begin{equation}
    e(S,S^{c}) \geq \frac{2(1-c)}{c}\cdot e
\end{equation}
To sum up, we have:
\begin{equation}
    e(S,S^{c}) \geq \begin{cases}
        \frac{2(1-c)}{1+c}\cdot min(e(S),e(S^{c})) + \frac{c}{4}|V|, &\mbox{if } min(e(S),e(S^{c}))> \frac{c^2}{4(1-c)}|V| \\
        \frac{2(1-c)}{c} \cdot min(e(S),e(S^{c})), &\mbox{if } min(e(S),e(S^{c}))\leq \frac{c^2}{4(1-c)}|V|
    \end{cases}
\end{equation}

\bibliographystyle{alpha}
\bibliography{sample}

\end{document}